\documentclass{amsart}
\pdfoutput=1
\usepackage{caption}
\usepackage{subcaption}
\usepackage{amssymb}
\usepackage{graphicx}
\usepackage{wrapfig}
\begin{document}

\author{G\'abor Fejes T\'oth}
\address{Alfr\'ed R\'enyi Institute of Mathematics,
Re\'altanoda u. 13-15., H-1053, Budapest, Hungary}
\email{gfejes@renyi.hu}

\author{L\'aszl\'o Fejes T\'oth}

\author{W{\l}odzimierz Kuperberg}
\address{Department of Mathematics \& Statistics, Auburn University, Auburn, AL36849-5310, USA}
\email{kuperwl@auburn.edu}

\title{Neighbors}
\thanks{The English translation of the book ``Lagerungen in der Ebene,
auf der Kugel und im Raum" by L\'aszl\'o Fejes T\'oth will be
published by Springer in the book series Grundlehren der
mathematischen Wissenschaften under the title
``Lagerungen---Arrangements in the Plane, on the Sphere and
in Space". Besides detailed notes to the original text the
English edition contains eight self-contained new chapters
surveying topics related to the subject of the book but not
contained in it. This is a preprint of one of the new chapters.}

\begin{abstract}
Two members of a packing are neighbors if they have a common boundary
point. A multitude of problems arises in connection with neighbors in a
packing. The oldest one concerns a dispute between Newton and Gregory about
the maximum number of neighbors a member can have in a packing of congruent
balls. Other problems ask for the average number of neighbors or the
maximum number of mutually neighboring members in a packing. The present work
gives a survey of these problems.
\end{abstract}

\maketitle

Two members of a packing are {\it{neighbors}} if they have a common boundary point.
In this chapter we survey results connected with the number of neighbors in a
packing.

\section{The Newton number of convex disks}

Let $N(S)$ denote the maximum of the number of neighbors that one member can have
in a packing of congruent copies of $S$. After {\sc{L.~Fejes T\'oth}}~\cite{FTL69b}
we call $N(S)$ the {\it{Newton number}} of $S$. The name racalls the dispute between
Newton and Gregory about the maximum number of congruent balls that can touch
another one of the same size without overlapping with each other. A rigorous
proof settling the dispute in favor of Newton was given by {\sc{Sch\"utte}} and
{\sc{van der Waerden}} \cite{Schuttevan der Waerden}. An alternative
name for $N(S)$ is the {\it{kissing number}} of $S$.

Some experiments with regular polygonal disks point to the conjecture that the
Newton number of a regular $n$-gon is: $12$ for $n=3$; $8$ for $n=4$; and $6$
for $n\ge5$. Indeed, with the exception of $n=5$, this was proved by
{\sc{B\"or\"oczky}}~\cite{Boroczky71} (see also {\sc{Youngs}}~\cite{Youngs} for the case of a
square, {\sc{Klamkin}}, {\sc{Lewis}} and {\sc{Liu}}~\cite{KlamkinLewisLiu} for the cases
$n=3,\ 4$ and $6$, and {\sc{Zhao}}~\cite{Zhao} for the case $n>6$). The Newton
number of the regular pentagon was proved to be six by {\sc{Linhart}}~\cite{Linhart73a} and
independently by {\sc{Pankov}} and {\sc{Dolmatov}} \cite{PankovDolmatov77,PankovDolmatov79}
and {\sc{Zhao}} and {\sc{Xu}}~\cite{ZhaoXu}.

{\sc{Schopp}}~\cite{Schopp70} proved that the Newton number of any disk of constant
width is at most $7$. The Newton number of the Reuleaux triangle is equal to
$7$ (Figure~1). {\sc{Kemnitz}} and {\sc{M\"oller}}~\cite{KemnitzMoller} determined
the Newton number of all rectangles. {\sc{L.~Fejes T\'oth}}~\cite{FTL67d} proved the
inequality
$$N(S)\le(4+2\pi)\frac{D}{w}+2+\frac{w}{D}$$
for a convex disk $S$ with diameter $D$ and width $w$. This estimate is exact in
many cases. For example, for the isosceles triangle $\Delta$ with
$\frac{w}{D}=\sin\frac{\pi}{19}$ we get $N(\Delta)\le 64$. But a simple construction
(Figure~2) shows that $N(\Delta)\ge 64$, hence $N(\Delta)=64$.

Other bounds for the Newton number of a convex disk involving different
parameters were given by {\sc{Hortob\'agyi}}~\cite{Hortobagyi72,Hortobagyi75,Hortobagyi76b} and
{\sc{Wegner}}~\cite{Wegner92}. Each of them considered a generalization of the Newton
number produced by counting the number of congruent copies of a convex disk
$K$ that can touch another convex disk $C$. Wegner's result
includes the computation of the Newton number of the
$30^{\circ}$-$30^{\circ}$-$120^{\circ}$ triangle, which, as conjectured, turns
out to be 21. Further problems about Newton numbers are treated in
{\sc{Harborth}}, {\sc{Koch}} and {\sc{Szab\'o}}~\cite{HarborthKochSzabo}, {\sc{Kemnitz}}, {\sc{M\"oller}} and {\sc{Wojzischke}}
~\cite{KemnitzMollerWojzischke}, {\sc{Kemnitz}} and {\sc{Szab\'o}}~\cite{KemnitzSzabo} and
{\sc{Kemnitz}}, {\sc{Szab\'o}} and {\sc{Ujv\'ary-Menyh\'art}}~\cite{KemnitzSzaboUjvary-Menyhart}.

\medskip
\centerline {\immediate\pdfximage height4cm
{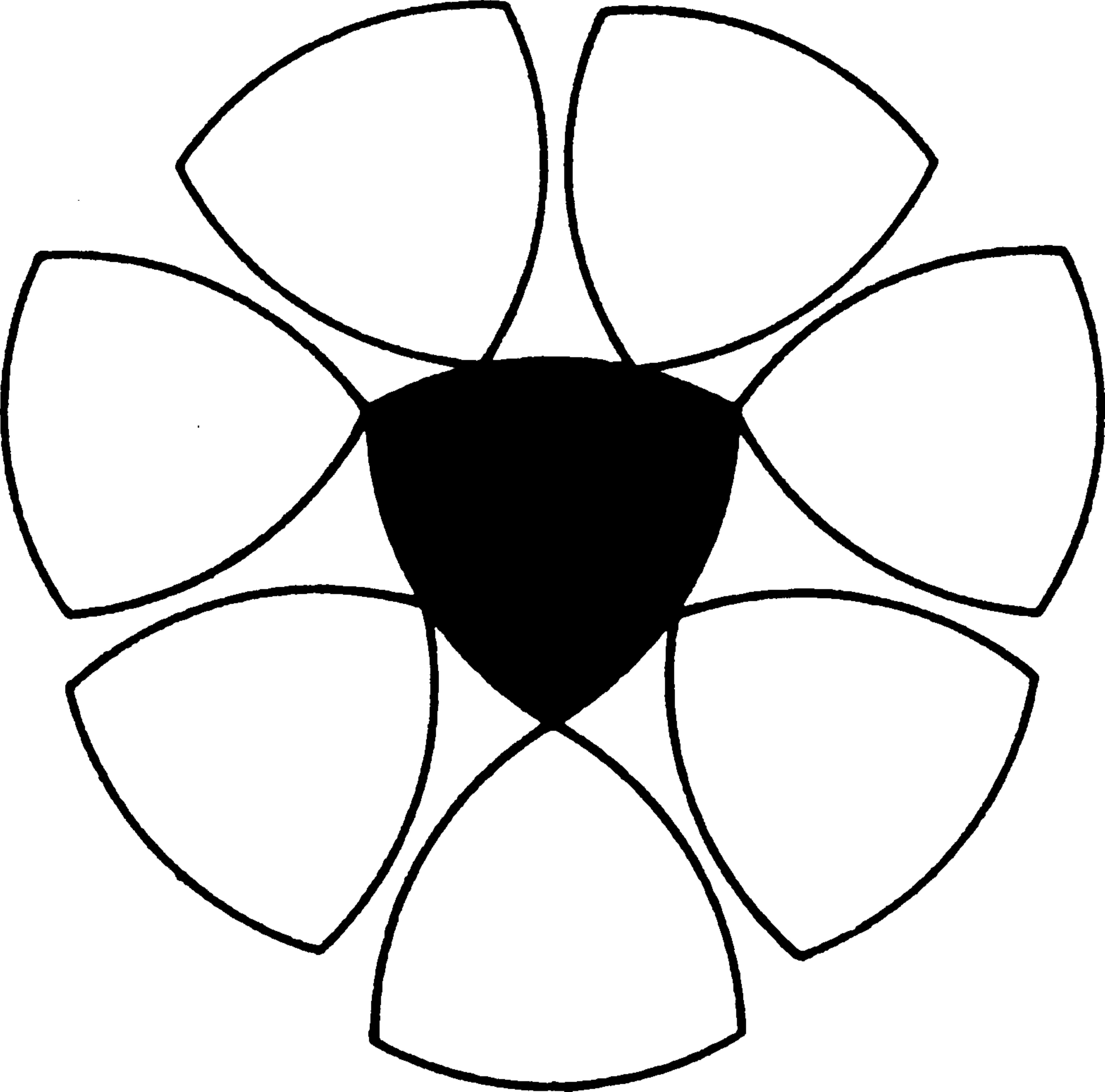}\pdfrefximage \pdflastximage\hskip2truecm
\immediate\pdfximage height7cm
{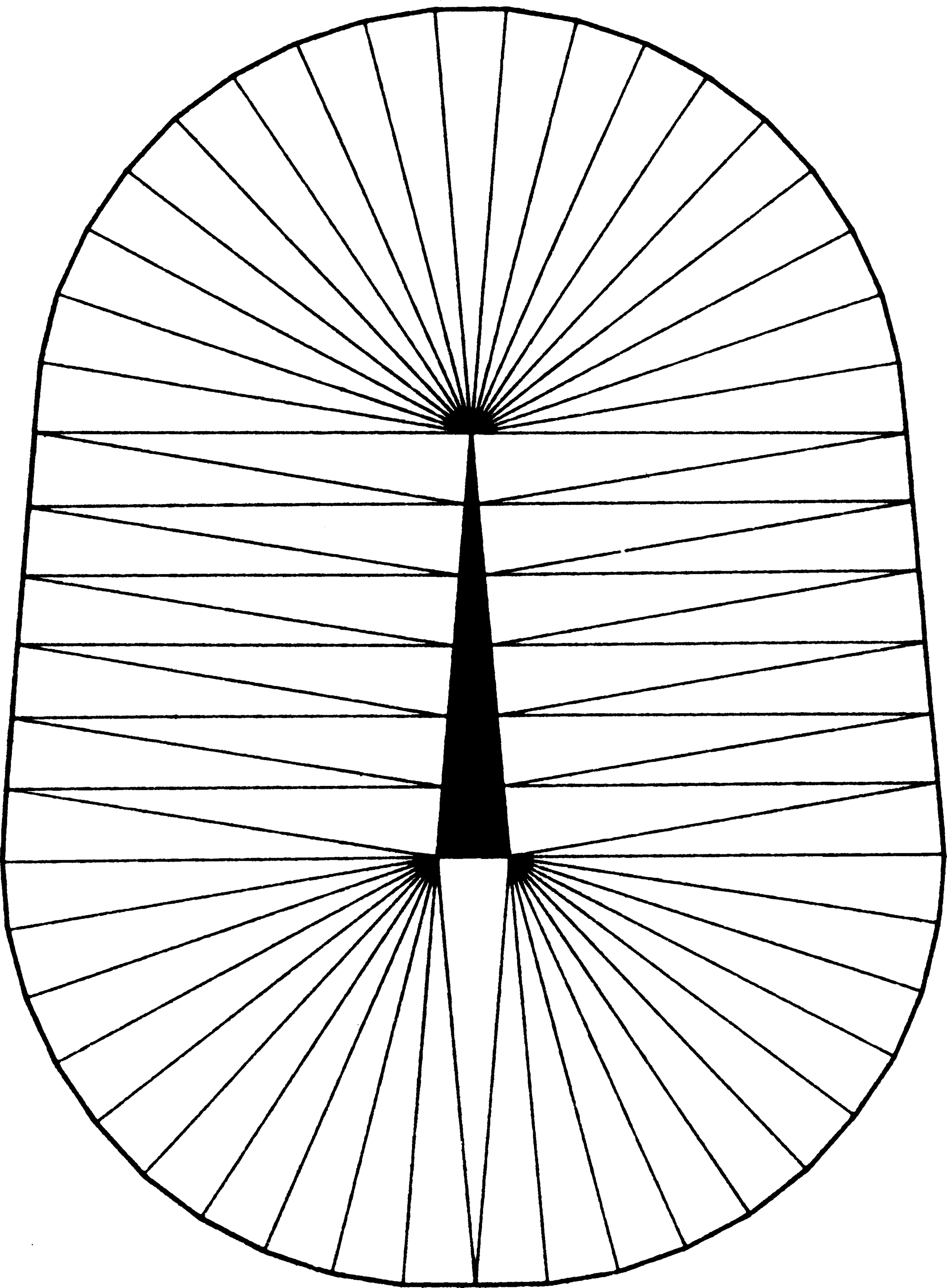}\pdfrefximage \pdflastximage}
\smallskip{\centerline{Figure~1 \hskip4.2truecm Figure~2}}

\section{The Hadwiger number of convex disks}

The {\it Hadwiger number} is defined in a similar way as the Newton number, but
considering only translated copies of $S$ instead of allowing all congruent ones.
The notion was so named by {\sc L.~Fejes T\'oth}~\cite{FTL70} because of the following
result of {\sc Hadwiger}~\cite{Hadwiger57c}: If each of $k$ mutually non-overlapping translates
of an $n$-dimensional convex body $K$ touches $K$, then $k\le3^n-1$.  The inequality
is sharp, as equality occurs for the parallelotope, and, as {\sc Groemer}~\cite{Groemer61a} proved,
only for the parallelotope (see also {\sc{Gr\"unbaum}}~\cite{Grunbaum61} for the planar case). The
Hadwiger number of a convex body $K$, denoted by $H(K)$, is alternately called the
{\it translative kissing number} of $K$.

In the plane, how large can the Hadwiger number be for a non-convex Jordan region?
While for each such region the number is finite, {\sc{Cheong}} and {\sc{Lee}}~\cite{CheongLee}
showed that, surprisingly, there exist Jordan regions with arbitrarily large Hadwiger number.
For starlike regions, however, {\sc{A.~Bezdek}}~\cite{BezdekA97} showed that $75$ is
an upper bound and {\sc{L\'angi}}~\cite{Langi11b} lowered
Bezdek's bound to 35. For the Hadwiger number of centrally symmetric starlike regions
{\sc{L\'angi}}~\cite{Langi09} established the stronger upper bound of $12$. It is still unknown
if the number can be greater than $8$. {\sc{L\'angi}}~\cite{Langi11b} determined the Hadwiger
number of a special class of regions. According to L\'angi a {\it{pocket}} of $R$ is a
connected component of ${\rm{conv}}\,R\setminus{R}$. In their examples, Cheong and Lee achieve
an arbitrarily large Hadwiger number only for a sequence of polygons whose number of pockets
is not bounded above. L\'angi determined the Hadwiger number of a region with one pocket;
depending on certain properties of the region, the number is either 6 or 8.

{\sc Boju} and {\sc Funar} \cite{BojuFunar93,BojuFunar11} investigated the
generalization of the Hadwiger number in which mutually non-overlapping,
$\lambda$-homothetic copies of $K$ for some $\lambda>0$ touch $K$ and they
gave corresponding bounds.

\section{Translates of a Jordan disk with a common point}

An interesting related question was addressed by {\sc{A.~Bezdek, K.~Kuperberg}}
and {\sc{W.~Kuperberg}}~\cite{BezdekAKuperbergKuperberg}: What is the maximum
number of non-overlapping translates of a Jordan disk that can have a common
point? They proved that the number is at most $4$, and they characterized the
disks for which $4$ non-overlapping translates can have a common point.
Earlier, {\sc{K.~Kuperberg}} and {\sc{W.~Kuperberg}}~\cite{KuperbergKKuperbergW}
proved the same bound for starlike disks. In contrast, in space the
corresponding number can be arbitrarily large even for starlike solids.

\section{The number of touching pairs in finite packings}

A problem of {\sc{Erd\H{o}s}} \cite{Erdos46} asks for the maximum number of
occurrences of the minimum distance between $n$ points in the plane. Erd\H{o}s
proved the upper bound $3n-6$ and pointed out that the example of the regular
triangular lattice shows that the minimum distance can occur $3n-c\sqrt{n}$
times. Unaware of the work of Erd\H{o}s, {\sc{Reuter}}~\cite{Reuter} restated
the problem as a conjecture in the language of circle packings: The maximum
number of touching pairs among $n$ unit circles forming a packing in the plane
is $\lfloor{3n-\sqrt{12n-3}}\rfloor$. The conjecture was verified by
{\sc{Harborth}}~\cite{Harborth}. {\sc{Kupitz}} \cite{Kupitz} gave a complete
description of the extremal packings; they are all subsets of the triangular
lattice. {\sc{Brass}}~\cite{Brass96} extended Harborth's
result to packings consisting of translates of a convex disk different from
a parallelogram and showed that the corresponding number for parallelograms
is $\lfloor{4n-\sqrt{28n-12}}\rfloor$. The papers by {\sc{ K. Bezdek}}
~\cite{BezdekK02,BezdekK12a} and {\sc{ K. Bezdek}} and {\sc{Reid}}
\cite{BezdekKReid} give bounds for the number of touching pairs in packings
of congruent balls. The latter paper investigates the number of mutually
touching triples and quadruples as well. {\sc{Bowen}} \cite{Bowen00} studied
contact numbers in circle packings in the hyperbolic plane.
{\sc{K. Bezdek, B.~Szalkai}} and {\sc{I.~Szalkai}} \cite{BezdekKSzalkaiSzalkai},
{\sc{K. Bezdek}} and {\sc{Nasz\'odi}} \cite{BezdekKNaszodi18b},
{\sc{K. Bezdek, Khan}} and {\sc{Oliwa}} \cite{BezdekKhanOliwa}, and
{\sc{Nasz\'odi}} and {\sc{Swanepoel}} \cite{NaszodiSwanepoel22} investigated
the number of touching pairs in totally separable packings. For a survey
on contact numbers in packings see {\sc{K. Bezdek}} and {\sc{Khan}}
\cite{BezdekKKhan18b}.

\section{$n$-neighbor packings}

If every member of a packing has exactly $n$, or at least $n$ neighbors, we call it
an $n$-{\it neighbor} or an $n^+$-{\it neighbor} packing, respectively. An easy
construction shows that there exists a zero-density 5-neighbor packing of the
plane with translates of a parallelogram. It turned out that this property
characterizes parallelograms (see {\sc{L.~Fejes T\'oth}}~\cite{FTL73a}). For a
convex disk $K$ other than a parallelogram, {\sc{Makai}}~\cite{Makai87} proved that
every $5^+$-neighbor packing with translates of $K$ is of density greater than
or equal to $3/7$. Equality can occur only if $K$ is a triangle. For centrally
symmetric convex disks the corresponding lower bound is $9/14$, attained only for
affine regular hexagons. In \cite{Makai87} {\sc{Makai}} only sketched the
proof for this last statement; details are given in {\sc{Makai}}~\cite{Makai18}.

Concerning $6^+$-neighbor packings, the following is known: {\sc{L.~Fejes T\'oth}}
~\cite{FTL73a} showed that the density of a $6^+$-neighbor packing with translates
of a convex disk is at least $1/2$, and {\sc{Makai}}~\cite{Makai87} obtained the
corresponding lower bound of $3/4$ for centrally symmetric disks. Each of
these bounds is sharp, and, again, the extreme values are produced only by
triangles and affine regular hexagons, respectively. {\sc{Chv\'atal}}~\cite{Chvatal75}
proved that the density of a $6^+$-neighbor packing with translated
parallelograms is at least $11/15$.

According to {\sc{L.~Fejes T\'oth}} \cite{FTL73a} the density of a $5^+$-neighbour
packing of congruent circular disks is at least $\sqrt3\pi/7$. Concerning $5^+$
neighbor packing of non-congruent circles {\sc{G.~Fejes T\'oth}} and
{\sc{L.~Fejes T\'oth}} \cite{FTGFTL91} proved that there is a constant
$h_0=0.53329\ldots$ such that any $5^+$ neighbor packing of circular disks whose
homogeneity exceeds $h_0$ has positive density. The constant $h_0$ is the unique
real root of the equation $8h^3+3h^2-2h-1$. It cannot be replaced by a smaller one:
A $5^+$-neighbor packing of circles of homogeneity $h_0$ and density 0
is shown in Figure~3.

\medskip
\centerline {\immediate\pdfximage height3cm
{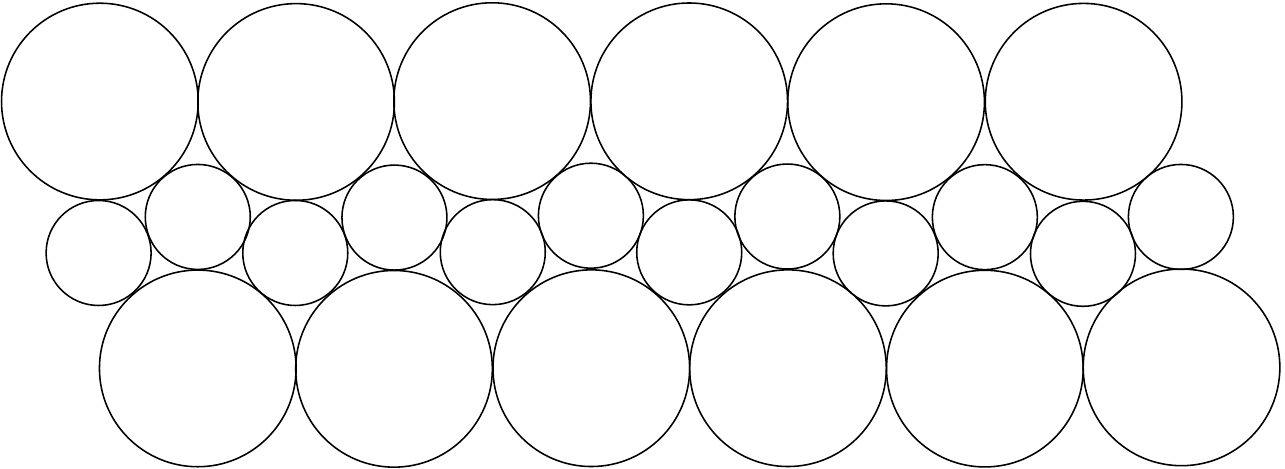}\pdfrefximage \pdflastximage}
\smallskip{\centerline{Figure~3}}
\medskip

{\sc{L.~Fejes T\'oth}}~\cite{FTL77} stated the following interesting conjecture:
The homogeneity of a $6^+$ neighbor packing of circles is
either 1 or 0. The conjecture was confirmed by {\sc{B\'ar\'any, F\"uredi}}
and {\sc{Pach}}~\cite{BaranyFurediPach}. They proved the somewhat stronger statement
that in a $6^+$-neighbor of circular disks either all disks are congruent or
arbitrarily small disks occur. Their proof combines a geometric idea with
a combinatorial one, each of interest on its own.

The angle at the boundary point $a$ of a convex disk $K$ is the
measure of the smallest angular region with apex $a$ containing $K$.
The minimum of the angle taken for all boundary points of $K$ is
called the \emph{minimal angle} of $K$. A further problem due to
{\sc{L.~Fejes T\'{o}th}}~\cite{FTL69c} asked for the numbers $n$
for which there is an $n$-neighbor packing consisting of convex
disks with given minimal angle $\pi/h$. {\sc{Linhart}}
\cite{Linhart74b} proved that for such a packing
$n\le\max\{5,2\lfloor{h}\rfloor-1\}$ and showed that this bound
is sharp for $\lfloor{h}\rfloor\le6$.

Let $t(K)$ denote the largest number $n$ for which there is a finite $n^+$-neighbor
packing of translates of a convex disk $K$, and let $m(K)$ denote the minimum cardinality
of such a packing. {\sc Talata}~\cite{Talata02} proved that $t(K)=4$ and
$m(K)=12$ if $K$ is a parallelogram, otherwise $t(K)=3$ and $m(K)=7$.

{\sc Wegner}~\cite{Wegner71} solved two problems posed by {\sc L.~Fejes T\'{o}th}~\cite{FTL69c}
concerning the existence of certain packings in which every member has the same number
of neighbors. On one hand, he constructed for every $n\ge3$ a stable $n$-neighbor packing of
congruent convex disks, and, on the other hand, he presented a $5$-neighbor packing consisting
of $32$ congruent smooth convex disks (see Figure~4). Independently from Wegner,
{\sc Linhart}~\cite{Linhart73b} also found a finite $5$-neighbor
packing of congruent smooth convex disks using a very similar idea.

\medskip
\centerline {\immediate\pdfximage height5.5cm
{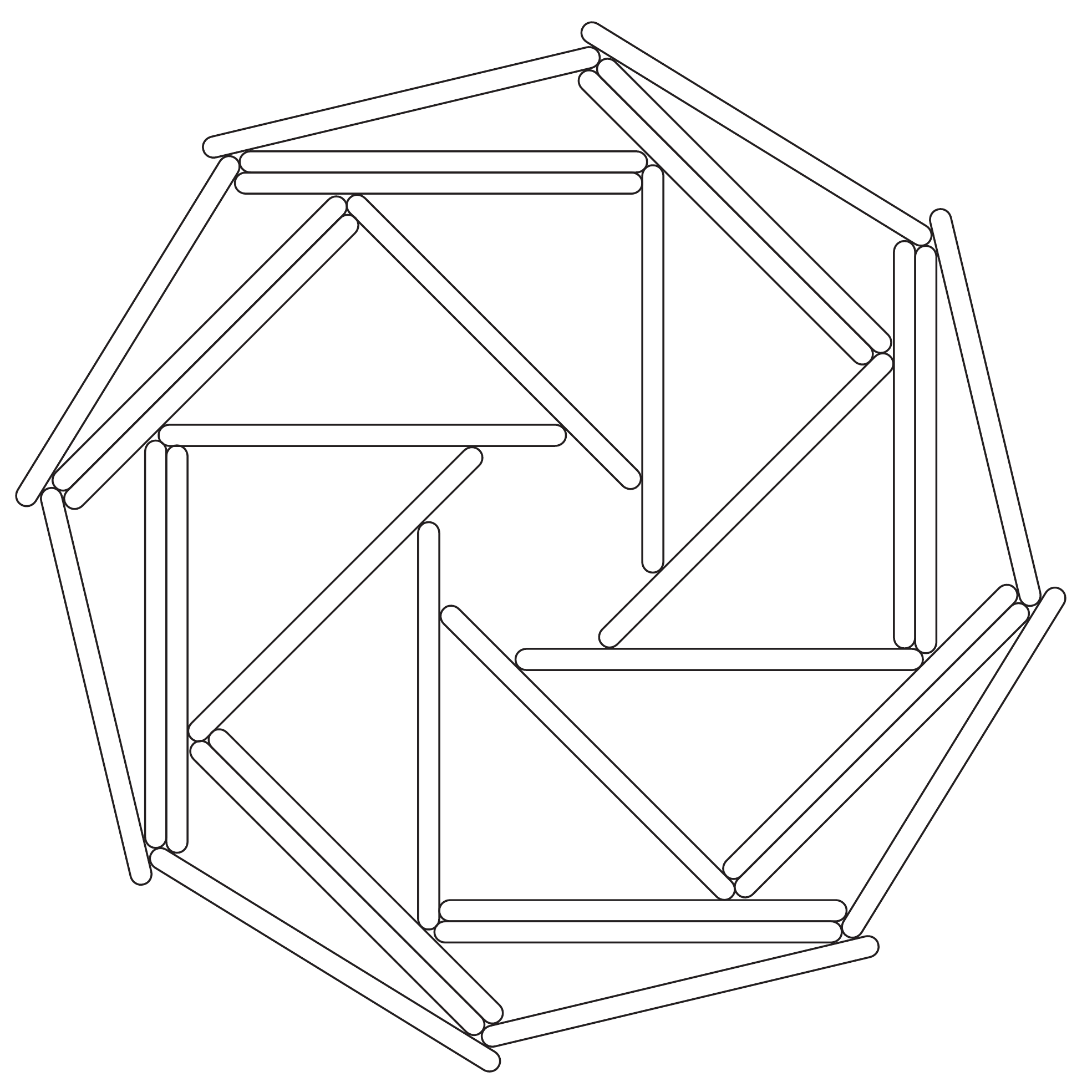}\pdfrefximage \pdflastximage}
\smallskip{\centerline{Figure~4}}

\section{Maximal packings}

A packing of congruent copies of a convex disks is said to be a {\it maximal
packing} if the number of neighbors of every disk equals its Newton
number. The faces of each of the tilings $\{6,3\}$, $\{4,4\}$ and $\{3,6\}$
form such a maximal packing.

{\sc{G\'acs}} \cite{Gacs} proved that there is an absolute constant $N$,
such that the number of neighbors in a maximal packing cannot exceed $N$.
In the tiling obtained by dissecting each face of $\{3,6\}$ into three
congruent triangles (Figure~5) every cell has $21$ neighbors. In view of
the result of {\sc{Wegner}} \cite{Wegner92} mentioned above, it follows that
this tiling is a maximal packing. On the other hand, {\sc{Linhart}}~\cite{Linhart77a}
showed that in no maximal packing of convex disks can the number of neighbors
exceed $21$.

\medskip
\centerline {\immediate\pdfximage width4cm
{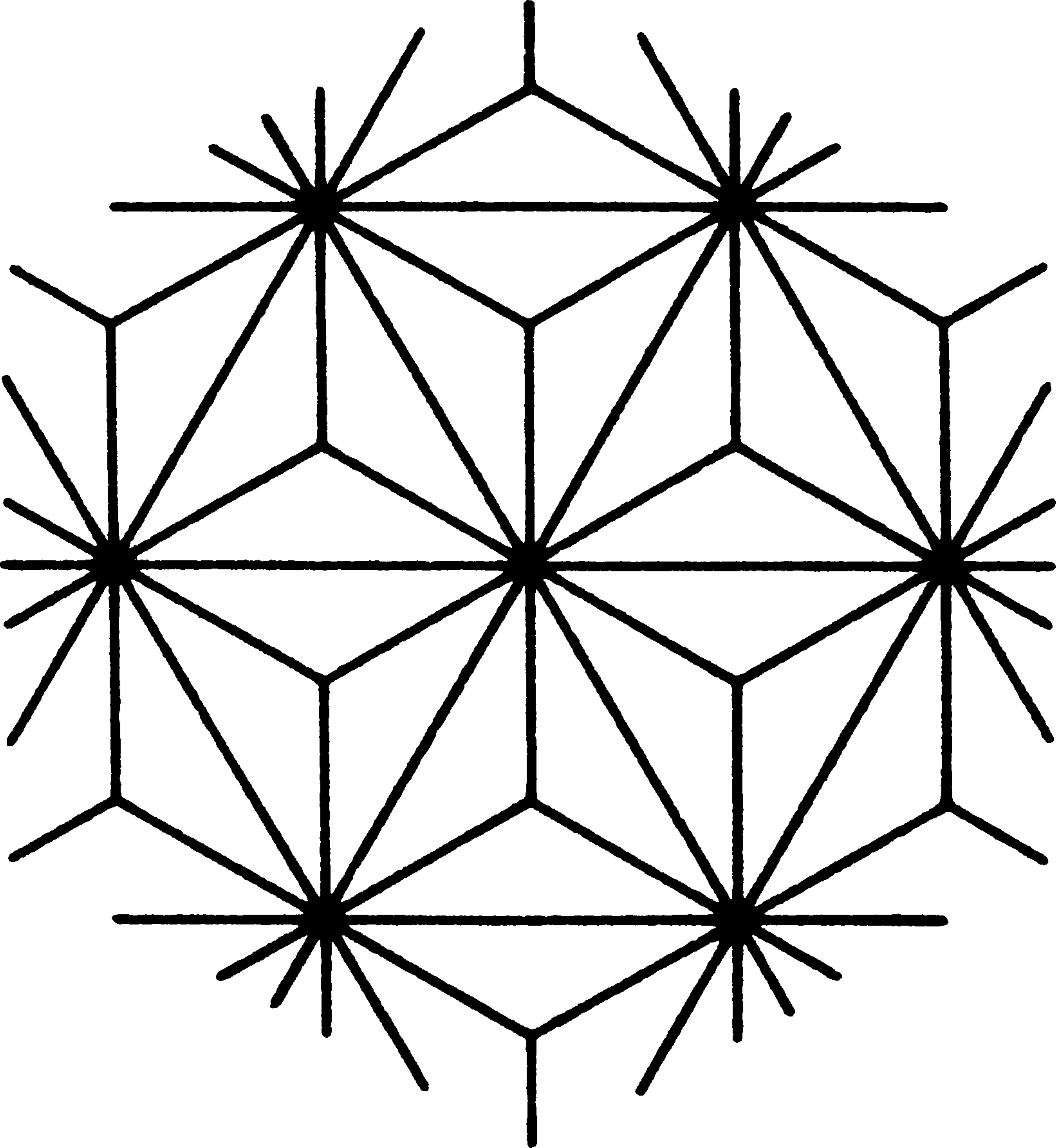}\pdfrefximage \pdflastximage}
\smallskip{\centerline{Figure~5}}
\medskip

{\sc{B\"or\"oczky}}~\cite{Boroczky71} proved that the faces of every regular
spherical tiling form a maximal packing. However at most $48$ regular
hyperbolic tilings are  maximal (see {\sc{Florian}}~\cite{Florian75b} and
{\sc{A.~Florian}} and {\sc{H.~Florian}}~\cite{FlorianFlorian75a,FlorianFlorian75b}).
{\sc{Linhart}}~\cite{Linhart75} showed that the tiling $\{7,3\}$ is maximal.

\section{Higher order neighbors}

We say that in a packing of convex disks each disk is the ``zero-th neighbor'' of itself. A disk
$B$ is the {\it$k$-th neighbor} of the disk $A$ if it is a neighbor of a $(k-1)$-th neighbor of $A$
other than a $j$-th neighbor of $A$ with $0\le{j}<k-1$. The {\it$k$-th Newton number} $N_k(D)$
of a convex disk $D$ is the maximum of the total number of its $j$-th neighbors for $j=1,\ldots,k$.
The densest lattice packing of circles yields the lower bound $N_k(B^2)\ge3k(k+1)$. One could
conjecture that, for small values of $k$, say for $k\le12$, $N_k(B^2)=3k(k+1)$ holds. This is
obvious for $k=1$ and was confirmed by {\sc{L.~Fejes T\'oth}} and {\sc{Heppes}}~\cite{FTLHeppes67} for $k=2$.
But {\sc{L.~Fejes T\'oth}}~\cite{FTL69a} showed that $N_{14}(B^2)\ge636$, which is greater than
$3\cdot14\cdot15=630$. He also determined the asymptotic behavior of $N_k(B^2)$ for large $k$,
namely
$$\lim_{k\to\infty}N_k(B^2)/k^2=2\pi/\sqrt3\,.$$

Higher order Hadwiger numbers are defined and treated in {\sc{L.~Fejes T\'oth}}
~\cite{FTL75a} and {\sc{L.~Fejes T\'oth}} and {\sc{Heppes}}~\cite{FTLHeppes77}.

\section{The Newton number of balls}

Alternates to Sch\"utte and Van der Waerden's proof of the equality $N(B^3)=12$
were given by {\sc{Leech}} \cite{Leech56}, {\sc{B\"or\"oczky}} \cite{Boroczky03},
{\sc{Anstreicher}} \cite{Anstreicher04}, {\sc{Musin}} \cite{Musin06a}, {\sc{Glazyrin}}
\cite{Glazyrin20b} and {\sc{Maehara}} \cite{Maehara07b}, the last one being perhaps
the most elementary among them.

{\sc{Flatley, Tarasov, Taylor}} and {{\sc{Theil}} {\cite{Flatley+} proved that the
maximum number of tangent pairs among twelve non-overlapping unit balls tangent
to a thirteenth unit ball is 24, attained only in the case when the centers of the
twelve balls are the vertices either of a cuboctahedron or of a twisted
cuboctahedron. A {\it{twisted cuboctahedron}} is obtained by cutting a cuboctahedron
into two parts by a plane containing 6 edges forming a regular hexagon, and
rotating one part by an angle of $\pi/3$ around the axis through the center of the
hexagon and perpendicular to its plane. {\sc{R. Kusner, W. Kusner, Lagarias}} and
{\sc{Shlosman}} \cite{KusnerKusnerLagariasShlosman} described the configuration
space of 12 non-overlapping equal spheres of radius $r$ touching a central unit
sphere. They also gave a nice survey of the history of the twelve spheres
problem and the Tammes problem.

The value of $N(B^4)$ was determined by {\sc Musin} \cite{Musin03,Musin08},
while $N(B^8)$ and $N(B^{24})$ were determined by {\sc Odlyzko} and {\sc Sloane}
\cite{OdlyzkoSloane} and, independently, by {\sc Leven\v{s}tein} \cite{Levenstein79}.
The cases $n=8$ and $n=24$ were resolved by means of the linear programming method,
and the case $n=4$ by its modification. The corresponding values are:
$$
N(B^4)=24,\quad N(B^8)=240,\quad {\rm and}\quad N(B^{24})=196560.
$$
Each of these Newton numbers is realized in the unique densest lattice
packing of balls in the corresponding dimension.  Moreover, as shown by {\sc
Bannai} and {\sc Sloane} \cite{BannaiSloane}, each of the arrangements of balls that
realize the Newton number for $n=8$ and $n=24$ is unique, which for $n=4$
is still only conjectured. The Kabatjanski\u{\i}--Leven\v{s}te\u{\i}n bound
yields $N(B^n)\le2^{0.4041d+o(d)}$, while the best known lower bound due to
{\sc{Jenssen}}, {\sc{Joos}} and {\sc{Perkins}} \cite{JenssenJoosPerkins} is $N(B^n)\ge(1+o(1))\sqrt{{3\pi}/{8}}\log({3\sqrt2}/{4})n^{3/2}(2\sqrt3/2)^n$.
{\sc{Dostert}} and {\sc{Kolpakov}} \cite{DostertKolpakov20}
gave upper bounds for the Newton number of balls in spherical and hyperbolic
space.

Let $N_k^n$ be the $k$-th Newton number of the $n$-dimensional ball. We have
$N^1_2=4$ and, as expected, $N^2_2=18$. {\sc{L.~Fejes T\'oth}} and {\sc{Heppes}}
\cite{FTLHeppes67} showed that $56\le N^3_2\le63$ and
$168\le N^4_2\le232$. The lower bound $N^3_2\ge56$ comes from the
enumeration of balls in the first and second neighborhood in the
trapezo-rhombododecahedral  packing. The outcome is $12+44=56$. Interestingly
enough, the corresponding outcome for the rhombododecahedral packing amounts
to just $12+42=54$. This could be of some significance to the still
unanswered question why certain metals form the trapezo-rhombododecahedral
structure in their crystals, while others form a rhombododecahedral structure.
In the densest lattice packing of the $4$-dimensional ball, found by
{\sc{Korkine}} and {\sc{Zolotareff}} \cite{KorkineZolotareff72}, the corresponding
number is $24+144=168$.

\section{$n$-neighbor packing of congruent balls}

What is the chromatic number of a finite packing of congruent balls? For which
numbers $n$ does there exist a finite $n$-neighbor packing of congruent balls?
Two consecutive layers from the densest lattice packing of balls is a
nine-neighbor packing of density 0. Does there exist a ten-neighbor packing of
zero density?  Motivated by these questions {\sc{L. Fejes T\'oth}} and
{\sc{Sachs}} \cite{FTLSachs} stated the following conjectures.

``Conjecture A. The maximum number of points which can be placed on an open unit
hemisphere with at least unit distance from one another is equal to eight.

Conjecture B. The maximum number of points which can be placed on a closed unit
hemisphere with at least unit distance from one another is equal to nine.

Conjecture AB. If nine points lie on a closed unit hemisphere with at least unit
distance from one another then six of them are on the boundary of the hemisphere.''

Since the condition that six points lie on the boundary uniquely determines the
configuration, Conjecture~AB implies the other two. It follows from Conjecture~A
that the chromatic number of a finite packing of unit balls is at most 9, while
the highest known chromatic number of a finite packing of unit balls is 5 (see
{\sc{Maehara}} \cite{Maehara07a}). It also follows that no finite nine-neighbor
packing of congruent balls exists. Conjecture~B implies that the density of a
ten-neighbor packing is positive.

Conjecture B was confirmed by {\sc{G. Fejes T\'oth}} \cite{FTG81}. Alternative
proofs based on the difficult result that the Newton number of the 3-ball is
12 were given by {\sc{Sachs}} \cite{Sachs} and  by {\sc{A.~Bezdek}} and
{\sc{K.~Bezdek}} \cite{BezdekABezdekK88}. Finally, {\sc{Kert\'esz}}
\cite{Kertesz94} proved Conjecture AB. The analogous problem was
investigated in higher dimensions as well. Let $B(n)$ be the maximum
number of points which can be placed on a closed unit hemisphere in $E^n$
with at least unit distance from one another. Equivalently, $B(n)$ is the
maximum number of non-overlapping unit balls that can touch another unit ball
at points of a closed hemisphere. Corresponding to this, $B(n)$ is called
the {\it{one-sided kissing number}} of $B^n$. {\sc{Szab\'o}} \cite{Szabo91}
proved that $B(4)\le20$. Based on an extension of the Delsarte method
{\sc{Musin}} \cite{Musin06b} lowered this bound to the sharp value 18.
{\sc{Bachoc}} and {\sc{Vallentin}} \cite{BachocVallentin09} proved that
$B(8)=183$ and {\sc{Dostert, de Laat and Moustrou}} \cite{DostertdeLaatMoustrou}
proved uniqueness of the arrangement. Clearly, a finite $k^+$-packing of congruent
balls in $E^n$ exists only for $k<B(n)$. It appears, however, that the maximum
number $k$ for which there exists a finite $k^+$-packing of congruent balls in $E^n$
is considerably less than $B(n)$. {\sc{Alon}} \cite{Alon} gave an explicit
construction of a finite $2^{\sqrt{n}}$-packing of unit balls.

Unit balls centered at the vertices of a triangle, tetrahedron, octahedron,
and dodecahedron form a two-, three-, four-, and five-neighbor packing,
respectively. It is easy to see that for $n\le4$ these are the arrangements of
minimal cardinality. {\sc{G. Fejes T\'oth}} and {\sc{Harborth}}
\cite{FTGHarborth} showed that a five-neighbor packing of congruent balls has
at least 12 members, and reported a six-neighbor packing consisting of 240
congruent balls constructed by Gerd Wegner (see Figure~6). It is an
open question whether a finite seven-neighbor or eight-neighbor packing of
congruent balls exists.

It is natural to ask for the minimum density of an $n$-neighbor packing of
congruent balls for $n=10,\ 11$ and $12$.  For $n=10$ and 11 our knowledge is
limited. Even the question whether an 11-neighbor packing of congruent balls
exists is open. The case $n=12$ is of special interest. {\sc L.~Fejes
T\'{o}th} \cite{FTL69b,FTL89} conjectured that if in a packing with congruent
balls each ball is touched by exactly twelve other balls, then the packing
must consist of parallel hexagonal layers. This long-standing conjecture was
verified by {\sc{Hales}} \cite{Hales12b,Hales13}. {\sc{B\"or\"oczky}} and
{\sc{Szab\'o}} \cite{BoroczkySzabo15} gave an alternative proof based on the
result of {\sc{Musin}} and {\sc{Tarasov}} \cite{MusinTarasov12} about the
densest packing of 13 spherical caps.

\medskip
\centerline {\immediate\pdfximage width6cm
{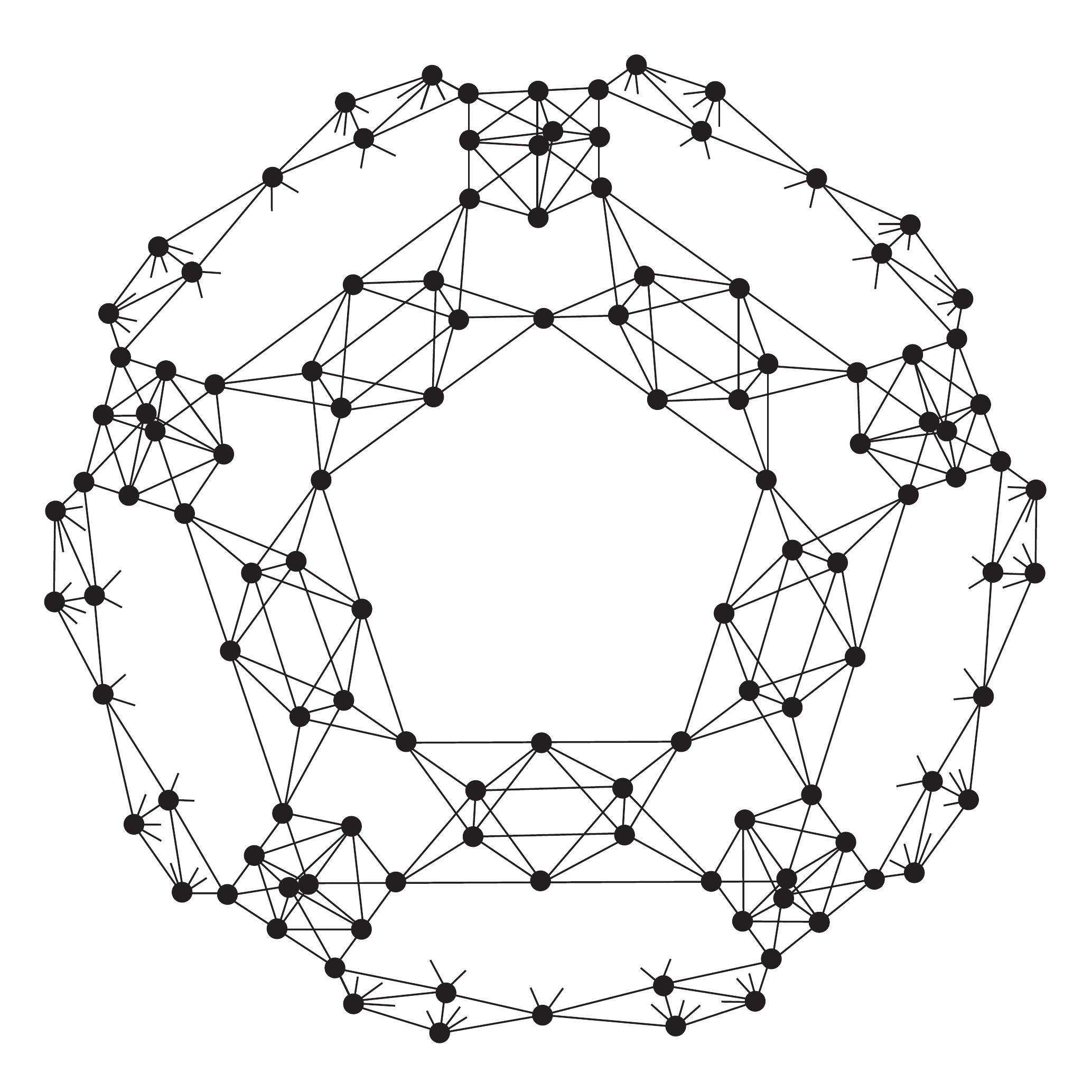}\pdfrefximage \pdflastximage}
\smallskip{\centerline{Figure~6}
\medskip

{\sc{Harborth, Szab\'o}} and {\sc{Ujv\'ary-Menyh\'art}}
\cite{HarborthSzaboUjvary-Menyhart} dropped the condition that the balls
are congruent, and constructed finite $n$-neighbor packings of balls for
all $n\le12$ except for 11. The question whether a finite 11-neighbor
packing of balls exists remains open.

Since the smallest ball in a finite ball-packing has at most twelve
neighbors, there is no finite $n$-neighbor packing for $n\ge12$. On the
other hand, the average number of neighbors in a finite packing of balls
can be greater than 12. {\sc{G.~Kuperberg}} and {\sc{Schramm}}
\cite{KuperbergSchramm} constructed a finite packing of balls in which
the average number of neighbors is $666/53= 12.566\ldots$. A packing of
balls with a slightly greater average number of neighbors of
$7656/607=12.612\ldots$ was given by {\sc{Eppstein, G. Kuperberg}}
and {\sc{Ziegler}} \cite{EppsteinKuperbergZiegler}.  {\sc{G.~Kuperberg}}
and {\sc{Schramm}} \cite{KuperbergSchramm} proved the upper bound
$8+4\sqrt3 = 14.928\ldots$ for the average number of neighbors in a
packing of balls, which was improved to $13.955$ by {\sc{Glazyrin}}
\cite{Glazyrin20a} and to $13.606$ by {\sc{Dostert, Kolpakov}} and
{\sc{Oliveira Filho}} \cite{DostertKolpakovOliveira}. The latter two
papers also give upper bounds for the average number of neighbors in
ball packings in higher dimensions.

{\sc{K.~Bezdek}}, {\sc{Connelly}} and {\sc{Kert\'esz}} \cite{BezdekKConnellyKertesz}
investigated packings of circles of radius $r$ on the sphere and proved that
there are positive numbers $\varepsilon$ and $r_0$ such that for $r\le r_0$
the average number of neighbors in such packings is at most $5-\varepsilon$.

\section{Results about convex bodies}

{\sc{Talata}} \cite{Talata02} investigated packings of translates of cylinders
in $E^3$ and showed that the maximum $n$ for which there exists a finite $n^+$-neighbor
packing of translates of a cylinder $C$ is 10 if $C$ is not a parallelepiped
and 13 if $C$ is a parallelepiped. He constructed a $10^+$-neighbor packing of
172 translates of a cylinder if the base is not a parallelogram and a
$13^+$-neighbor packing of 382 translates of a parallelepiped.

{\sc{L. Fejes T\'oth}} and {\sc{Sauer}} \cite{FTLSauer} proved that if in
a packing of translates of an $n$-dimensional cube, for each cube the total
number of its $j$-th neighbors for $0\le{j}\le{k}$ is more than
$(k+1)(2k+1)^{n-1}$, then the packing has positive density. The result is
sharp, as shown by the example of two consecutive layers in the grid of
cubes in which every cube has a total number of $(k+1)(2k+1)^{n-1}$
$j$-th neighbors with $0\le{j}\le{k}$. {\sc{K. Bezdek}} and {\sc{Brass}}
\cite{BezdekKBrass} generalized the case $k=1$ of this result for packings
by translates of an arbitrary $n$-dimensional convex body.

{\sc Zong} \cite{Zong96} constructed a lattice packing of tetrahedra in which
every tetrahedron touches $18$ others, conjectured that for a tetrahedron $T$
$H(T)\le18$, and proved that $H(T)\le19$. Since in a lattice packing each
member has an even number of neighbors, it follows that the number of neighbors
of a member in a lattice packing of tetrahedra cannot exceed 18. {\sc Talata}
\cite{Talata98b} gave a simple alternative proof for this. Later, {\sc Talata}
\cite{Talata99a} succeeded in proving Zong's conjecture about the Hadwiger
number of tetrahedra. Moreover, he showed that the packing of 18 translates
of a tetrahedron touching a nineteenth one is unique. He also gave a description
of all possible packings of 17 translates of a tetrahedron touching an
eighteenth one. He applied this result for the determination of the minimum and
maximum densities of $17^+$-neighbor translative packings of tetrahedra.
The Hadwiger number of the octahedron is 18. This was proved independently by
{\sc{Robins}} and {\sc{Salowe}} \cite{RobinsSalowe}, {\sc{Talata}}
\cite{Talata99c} and {\sc{Larman}} and {\sc{Zong}} \cite{LarmanZong}. The latter
author showed that the Hadwiger number of the rhombic dodecahedron is also 18.

{\sc Gr\"{u}nbaum} \cite{Grunbaum61} proved that the Hadwiger number of a
convex disk is always attained in a lattice packing. {\sc Zong} \cite{Zong94}
showed that this statement does not hold in any dimension $n\ge3$, namely it
fails for a cube truncated at some of its vertices. In the same article,
Gr\"{u}nbaum conjectured that the Hadwiger number of every convex body is
always even. Disproving the conjecture, {\sc Jo\'{o}s} \cite{Joos08a}
constructed a $3$-dimensional convex body whose Hadwiger number is $15$.

{\sc Talata} \cite{Talata98a} proved that there is an absolute constant $c>0$
such that $H(K)\ge 2^{cn}$ for every $n$-dimensional convex body $K$.
In \cite{Talata00b} {\sc Talata} gave the bound  $H(S_n)\ge
1.13488^{(1-o(1))n}$ for the $n$-dimensional simplex $S_n$, and for strictly
convex bodies $K$ he gave in \cite{Talata05} the following explicit bound:
$H(K)\ge\frac{16}{35}7^{(n-1)/2}$. {\sc{Robins}} and {\sc{Salowe}}
\cite{RobinsSalowe}, {\sc{Swanepoel}} \cite{Swanepoel99}, {\sc{Larman}} and
{\sc{Zong}} \cite{LarmanZong} and {\sc{Xu}} \cite{Xu} gave lower bounds for
the Hadwiger number of superballs.

{\sc Zong} \cite{Zong97a} studied the Hadwiger number of Cartesian products of
convex bodies $K$ and $L$, and he proved that if ${\rm dim}L\le2$, then
$H(K\times L)=(H(K)+1)(H(L)+1)-1$.  {\sc Talata} \cite{Talata05} showed that
if the dimension $L$ is higher than $2$ the equality does not always hold.

It is clear that the density of a saturated packing of congruent copies of
a convex body with a large average number of neighbors cannot be arbitrarily
small. {\sc{Groemer}} \cite{Groemer61c} proved that the density of a saturated
packing of translates of a convex body in $E^n$ with average number of
neighbors $\mu$ is at least $\frac{1}{3^n-\mu}$.

\section{Mutually touching translates of a convex body}

A problem of {\sc{Erd\H{o}s}} \cite{Erdos48} asked for the greatest cardinality
of a set of points in $E^n$ with the property that no angle determined by
three points is greater than $90^\circ$. Another problem posed by {\sc{Klee}}
asked for the largest antipodal set in $E^n$. Two points of a set are
{\it{antipodal}} if there are two parallel supporting hyperplanes of the
set, each containing one of the points while the whole set is contained
in the closed slab bounded by the supporting planes. A set is {\it{antipodal}}
if every pair of its points are antipodal.

{\sc{Danzer}} and {\sc{Gr\"unbaum}} \cite{DanzerGrunbaum} proved that the two
problems are equivalent and confirmed the conjecture of Erd\H{o}s that the
answer is $2^n$. Moreover, they showed that these problems are also equivalent
to the problem of finding the largest family of mutually touching
translates of a convex body in $E^n$. The {\it{touching number}} $t(K)$ of
a convex body $K$ in $E^n$ is the maximum number of pairwise touching
translates of $K$. Thus, the theorem of Danzer and Gr\"unbaum states that
$t(K)\le2^n$ with equality  attained only for a parallelotope.

K\'aroly Bezdek and J\'anos Pach (see \cite[p.~98, Conjecture 13]{BrassMoserPach})
conjectured that if $C$ is a centrally symmetric convex body, then even a family
of pairwise touching homothetic copies of $C$ has at most $2^n$ members.
{\sc{Nasz\'{o}di}} \cite{Naszodi06} proved the upper bound $2^{n+1}$ without
assuming symmetry of $C$ and replacing ``homothetic'' with ``positively-homothetic''.
For centrally symmetric bodies Nasz\'{o}di's bound was lowered by {\sc{L\'{a}ngi}} and
{\sc Nasz\'{o}di} \cite{LangiNaszodi09} to $3\times2^{n-1}$. {\sc{F\"oldv\'ari}}
\cite{Foldvari} proved that the maximum number of pairwise touching positive homothetic
copies of a convex disk is $4$.

There is a fourth equivalent formulation of the problem of Erd\H{o}s and Klee. A
subset $S$ of a metric space $M$ is said to be {\it{equilateral}} provided each
pair of points in $S$ have the same distance. The maximum number of elements
in an equilateral set in $M$ is denoted by $e(M)$. If $K$ is the unit ball of
a Minkowski space $M$, then the set of mutually touching translates of $K$
corresponds to a set of equilateral points in $M$, and vice versa. Thus, for
an $n$-dimensional Minkowski space $M$, $e(M)\le2^n$, as was noted by
{\sc{Petty}} \cite{Petty71} and {\sc{P. S. Soltan}} \cite{SoltanPS}.

For $1\le{p}\le\infty$ and $n\ge1$ let $l_p^n$ and $l_p^\infty$ denote $R^n$
endowed with the norm $\|(x_1,\ldots,x_n\|=\left(\sum_{i=1}^n|x_i|^p\right)^{1/p}$
and $\|(x_1,\ldots,x_n\|=\max_{1\le{i}\le{n}}|x_i|$, respectively. In $l_1^n$
the standard basis vectors and their negatives form an equilateral set of
$2n$ points, and the set of standard basis vectors together with an appropriate
multiple of the all 1 vectors shows that $e(l^n_p)\ge{n+1}$. R.~B.~Kusner
(see {\sc{Guy}} \cite{Guy83}) conjectured that both examples are extremal, that is
$e(l^n_1)=2n$ and $e(l^n_p)={n+1}$ for $1<{p}<\infty$. The conjecture
concerning $l^n_1$ was confirmed for $n=3$ by {\sc{Bandelt, Chepoi}} and
{\sc{Laurent}} \cite{BandeltChepoiLaurent} and for $n=4$ by {\sc{Koolen, Laurent}}
and {\sc{Schrijver}} \cite{KoolenLaurentSchrijver}. Thus, the touching number of
an octahedron is 6, and the touching number of a cross-polytope in $E^4$ is 8.
Besides the Euclidean case $l^n_2$ settled by {\sc{C. Smith}} \cite{SmithC01},
Kusner's conjecture was confirmed for $l^n_4$ by {\sc{Swanepoel}} \cite{Swanepoel04a}
who also proved that the conjecture is false for all $1<p<2$ and sufficiently
large $n$, depending on $p$. However, it follows by continuity that, for fixed
$n$, if $p$ is close to 2 or 4 then $e(l^n_2)=n+1$ and $e(l^n_4)=n+1$.
{\sc{C. Smith}} \cite{SmithC01} and {\sc{Swanepoel}} \cite{Swanepoel14} gave explicit
bounds for $p$: If $|p-2|<\frac{2\log(1+2/n)}{\log(n+2)}$ then $e(l^n_2)=n+1$ and
if $|p-4|<\frac{4\log(1+2/n)}{\log(n+2)}$ then $e(l^n_4)=n+1$.
{\sc{C. Smith}} \cite{SmithC13} proved that there exists a constant $c_p$ such that
$e(l_p^n)\le c_pn^{(p+1)/(p-1)}$. Extending Smith's method, {\sc{Alon}} and
{\sc{Pudl\'ak}} \cite{AlonPudlak} proved $e(l_p^n)\le c_pn^{(2p+1)/(2p-1)}$,
and for odd integers $p\ge1$ established the bound $e(l_p^n)\le c_pn\log{n}$.

Petty conjectured that $e(M)\ge{n+1}$ for every $n$ dimensional Minkowski
space. In other terms the conjecture states, on one hand, that in
every $n$-dimensional Minkowski space there exists a full-dimensional
regular simplex and, on the other hand, that every convex body $K$ in $E^n$
admits a packing of $n+1$ mutually touching translates of $K$. Petty proved
the bound $e(M)\ge\min\{4,{n+1}\}$, confirming the conjecture for $n=3$.
Alternative proofs for the 3-dimensional case were given by {\sc{Kobos}}
\cite{Kobos13} and {\sc{V\"ais\"al\"a}} \cite{Vaisala}. This case is well
understood: {\sc{Gr\"unbaum}} \cite{Grunbaum63a} proved that $e(M)\le5$ if
the unit ball of $M$ is strictly convex, and {\sc{Sch\"urmann}} and
{\sc{Swanepoel}} \cite{SchurmannSwanepoel} proved that $e(M)\le6$ if the
unit ball is smooth. The latter authors gave an example of a smooth space
$M$ with $e(M)=6$ and also characterized the 3-dimensional Minkowski spaces
that admit equilateral sets of 6 and 7 points. The case $n=4$ of Petty's
conjecture was confirmed by {\sc{Makeev}} \cite{Makeev05}. For $n\ge5$ the
conjecture is still open, except for special classes of spaces.

The {\it{Banach-Mazur distance}} between two $n$-dimensional Minkowski
spaces is defined as $d(X,Y)=\inf\|T\|\|T^{-1}\|$, where the infimum is taken
over all linear, invertible operators $T$ from $X$ to $Y$. Let $M$ be an
$n$-dimensional Minkowski space. {\sc{Brass}} \cite{Brass90} and
{\sc{Dekster}} \cite{Dekster00a} proved that if $d(M,l_2^n)\le1+1/n$,
then $e(M)\ge n+1$. {\sc{Swanepoel}} and {\sc{Villa}}
\cite{SwanepoelVilla} verified Petty's conjecture also for the case that
$d(M,l_\infty^n)\le3/2$, and {\sc{Averkov}} \cite{Averkov} proved that even
$d(M,l_\infty^n)\le2$ guarantees $e(M)\ge n+1$. From these results
{\sc{Swanepoel}} and {\sc{Villa}} \cite{SwanepoelVilla} derived the lower
bound $e(M)\ge{e}^{c\sqrt{\log{n}}}$ using the theorem of {\sc{Alon}} and
{\sc{Milman}} \cite{AlonMilman} stating that for every $\varepsilon>0$
there exists a constant $c(\varepsilon)$ such that any $n$-dimensional
Minkowski space contains a subspace of dimension at least
$e^{c(\varepsilon)\sqrt{\log{n}}}$ whose Banach-Mazur distance to either
$l_2^n$ or $l_\infty^n$ is at most $\varepsilon$.

{\sc{Gonz\'alez Merino}} \cite{GonzalezMerino} verified Petty's conjecture
for spaces whose unit ball satisfies certain intersection properties and
{\sc{Kobos}} \cite{Kobos14} proved it if the unit ball is symmetric in each
of the hyperplanes $x_i=x_j$. Kobos also proved that the conjecture holds for
any $n-1$-dimensional subspace of $l_\infty^n$. {\sc{Frankl}} \cite{Frankl}
gave the following extension of this result: To every integer $k\ge2$
there is a bound $N(k)$ such that for $n>N(k)$ any $k$-dimensional subspace
of $l_\infty^n$ contains a set of $k+1$ equidistant points. As a corollary
she obtained that if a centrally symmetric polytope in $E^n$ has at most
$\frac{4}{3}n-\frac{1+\sqrt{8n+9}}{6}$ opposite pairs of facets, then
there are $n+1$ mutually translates of it. Unfortunately, the difference
body of the simplex has more faces, so this does not give a lower bound
for the touching number of simplices. {\sc{Koolen, Laurent}} and
{\sc{Schrijver}} \cite{KoolenLaurentSchrijver} gave examples of $n+2$
mutually touching $n$-dimensional simplices. Moreover, {\sc{Lemmens}}
and {\sc{Parsons}} \cite{LemmensParsons} proved that for $n\ge5$ and
$n\equiv1$ (mod $4$) the touching number of the $n$-dimensional simplex
is at least $n+3$.

Further reading about equilateral sets can be found in the papers by
{\sc{Swanepoel}} \cite{Swanepoel04b,Swanepoel18}.

\section{Mutually touching cylinders}

{\sc{Littlewood}} (\cite[Problem 7 on p.~20]{Littlewood}) asked the
following question: ˙˙Is it possible in 3-space for seven infinite circular
cylinders of unit radius each to touch all the others?" There are several
examples of six cylinders mutually touching each other. All known examples are
flexible with one degree of freedom. {\sc{Bozoki, Lee}} and
{\sc{R\'onyai}} \cite{BozokiLeeRonyai} answered Littlewood's question in the
affirmative. Fixing the angle of two cylinders the position of the remaining
five cylinders can be described by 20 parameters satisfying 20 multivariate
polynomial equations. Fixing the position of the initial cylinders
perpendicularly two essentially different approximate solutions were found by
numerical methods. Having found the approximate solutions, it was proved by
Smale's $\alpha$-theory (see {\sc{Smale}} \cite{Smale}), as well as by
interval-arithmetic computations, that the system of equations does indeed have real solution in
the neighborhood of the approximate solutions. Further numerical investigation
indicated that by fixing the angle of the initial two cylinders between any
angle $0<\phi<\pi$ gives a one-parametric class of solutions (see the
demonstration by {\sc{Scherer}} \cite{Scherer}). The two arrangements
found by fixing the angle perpendicularly are special elements of this
class. Figure~7 shows seven mutually touching cylinders.

Unaware of Littlewood's problem {\sc{Pikhitsa}} \cite{Pikhitsa} studied the
same question motivated by applications in physics. He also described a
configuration of seven mutually touching infinite cylinders numerically, without
giving a mathematically rigorous proof for its existence. {\sc{Pikhitsa}} and
{\sc{Choi}} \cite{PikhitsaChoi} gave numerical evidence for the existence of
nine mutually touching incongruent infinite cylinders.

W. Kuperberg presented a seemingly convincing physical model of eight
mutually touching congruent cylinders and asked whether the
cylinders really are mutually touching? {\sc{Ambrus}} and
{\sc{A. Bezdek}} \cite{AmbrusBezdek} showed that in this model there are two
cylinders that do not touch. The question whether there are eight mutually
touching congruent cylinders remains unanswered.

\medskip
\centerline {\immediate\pdfximage width8cm
{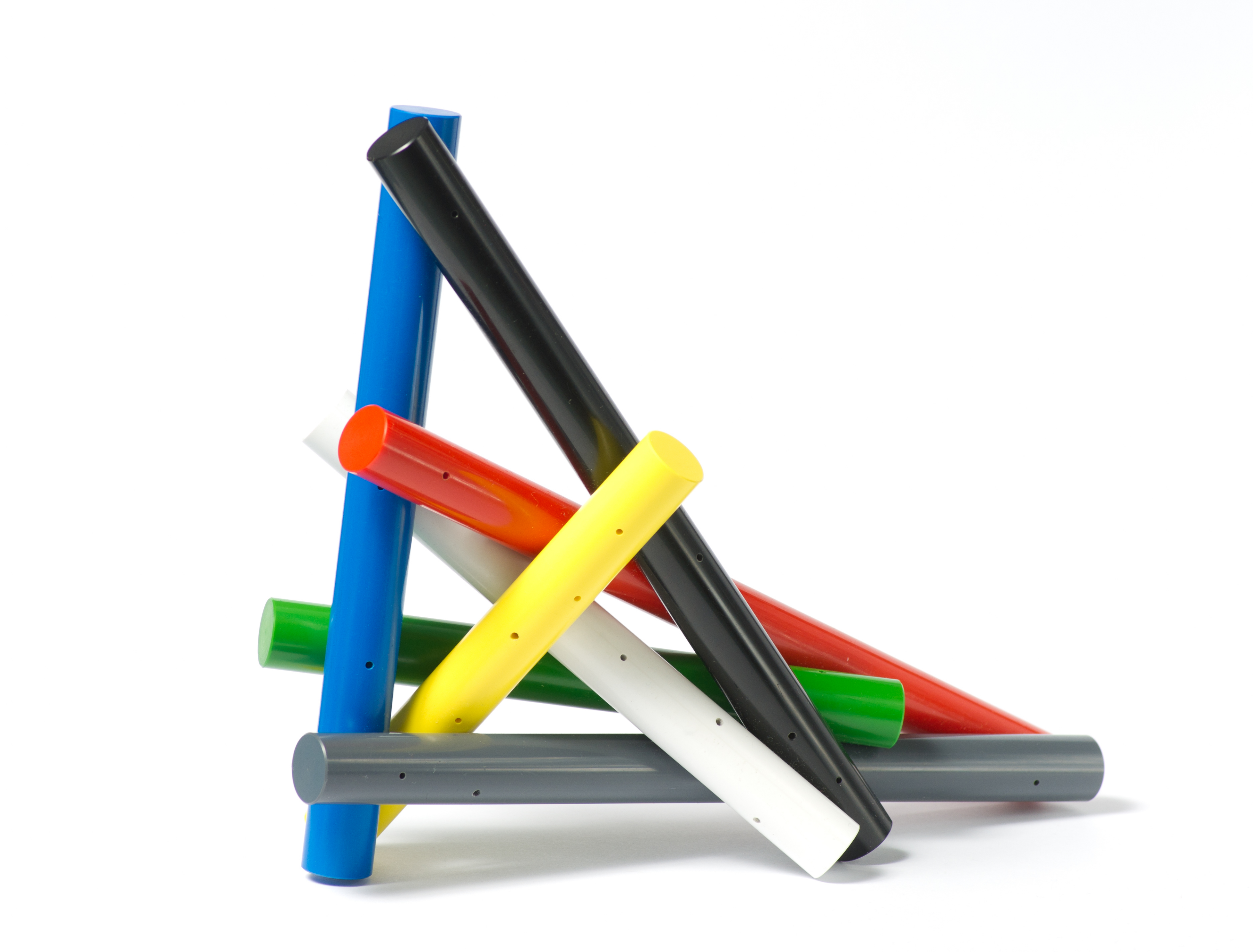}\pdfrefximage \pdflastximage}
\smallskip{\centerline{Figure~7}
\medskip

{\sc{A.~Bezdek}} \cite{BezdekA05a,BezdekA05b} gave two different proofs of the
statement that the number of mutually touching congruent infinite cylinders is
bounded. The first proof uses Ramsey's theorem, and the argument given there
shows that the number of cylinders is bounded even if we allow incongruent
cylinders. In the second paper he proved that no more than 24 congruent
infinite cylinders can mutually touch. {\sc{Pikhitsa}} and {\sc{Pikhitsa}}
\cite{PikhitsaPikhitsa17} also proved that the number of mutually touching infinite
cylinders of arbitrary base is bounded. Moreover, in \cite{PikhitsaPikhitsa19}
they claimed that no more than 10 infinite cylinders can mutually touch, and gave
numerical evidence for the existence of 10 mutually touching elliptic cylinders.

\section{Cylinders touching a ball}

{\sc{W.~Kuperberg}} \cite{Kuperberg90,Kuperberg14} asked for the maximum number of
unit-radius infinite cylinders touching a unit-radius ball. He conjectured
that the number in question is six, which can be realized in several different
ways. {\sc{Heppes}} and {\sc{Szab\'o}} \cite{HeppesSzabo} gave two different
proofs of the upper bound 8 on the number of cylinders. They also discussed the
same problem for higher dimensions, and for other radii of the touching
cylinders. {\sc{Brass}} and {\sc{Wenk}} \cite{BrassWenk} computed the
portion of area cut out by a cylinder touching the unit ball from a concentric
sphere of radius $\sqrt{4.7}$, which came out to be greater than $1/8$, showing
that the number of touching cylinders is at most 7. While the question whether 7
mutually disjoint infinite cylinders of unit radius can touch the unit ball
remains open, it turned out that 6 cylinders of radius $r>1$ can touch it. The
first such example with radius $r=1.049659$ was given by {\sc{Firsching}}
\cite{Firsching15} (see also {\sc{Firsching}} \cite{Firsching14}) by a
numerical exploration of the corresponding 18-dimensional configuration
manifold. {\sc{Ogievetsky}} and {\sc{Shlosman}} \cite{OgievetskyShlosman19a,
OgievetskyShlosman19b,OgievetskyShlosman19c,OgievetskyShlosman21a,
OgievetskyShlosman21b} devoted a series of papers to the study of the
configuration space of cylinders touching a ball. They found a packing of 6
infinite cylinders of radius $r=\frac{1}{8}(3+\sqrt{33})\approx1.093070331$
touching the unit ball, and believe that this value of the radius is the maximum.

{\sc{Starostin}} \cite{Starostin} investigated tubes touching a ball or
another tube.

\section{Neighbors in lattice packings}

Concerning the number of neighbors in densest lattice packings {\sc{Swinnerton-Dyer}}
\cite{Swinnerton-Dyer} proved that for every convex body $K$ in $E^n$ there is a
lattice packing of $K$ in which every member touches at least $n(n+1)$ others.
{\sc{M. J. Smith}} \cite{SmithMJ75} extended the result of Swinnerton-Dyer to
compact sets $S$ for which $S-S$ has non-empty interior. On the other hand,
{\sc Gruber} \cite{Gruber86} proved that, in the sense of Baire category,
typical convex bodies have at most $2n^2$ neighbors in their densest
lattice packings. The ball is not typical: {\sc{Vl\u{a}du\c{t}}}
\cite{Vladut19} constructed for a sequence of dimensions $n_i$ lattice
ball packings in $E^{n_i}$ in which the balls have $2^{0.0338n_i+o(n_i)}$
neighbors, and in \cite{Vladut21} he constructed sequences of lattice
packings of superballs in $E^n$ with an exponential number of neighbors.
The difference between the Hadwiger number and the maximum number $H_L(K)$
of neighbors of $K$ in a lattice packing of $K$ can be large. {\sc{Talata}}
\cite{Talata98b} proved that for every $n\ge3$ there exists an $n$-dimensional
convex body $K$ such that $H(K)-H_L(K)\ge2^{n-1}$.

{\sc{Groemer}} \cite{Groemer68b} studied the number of neighbors in connected
lattice packings and proved that in a thinnest connected lattice packing of
an $n$-dimensional convex body, each body has at least $2n$ and at
most $2(2^n - 1)$ neighbors.
\vfil\eject

\small{
\bibliography{pack}}